\documentclass[12pt]{iopart}

\usepackage{aas_macros}
\usepackage{harvard} 
\usepackage{graphicx}
\usepackage{fullpage}

\newcommand{\amc}{{(1-\mu) }}

\usepackage{amssymb}
\begin{document}

\title[Linearization .... Chermnykh's problem]{Linearization of the Hamiltonian around the triangular equilibrium points  in  the generalized photogravitational Chermnykh's problem}
\author{Badam Singh Kushvah}
\address{Department of    Mathematics, National Institute of Technology,  Raipur (C.G.) -492010, INDIA}
\ead{bskush@gmail.com}
\begin{abstract}
 Linearization of the Hamiltonian is being performed around the triangular equilibrium points  in  the generalized photogravitational Chermnykh's problem. The  bigger  primary  is being considered as a source of radiation and small primary as an oblate spheroid. We have found the normal form of the second order part of the Hamiltonian. For this we have solved the aforesaid set of equations. . The effect of radiation pressure,  gravitational potential from the belt on the linear stability have been  examined analytically and numerically. 
\end{abstract}

\submitto{\JPD}
\maketitle
\section{Introduction}
\label{intro}
The Chermnykh's problem is new kind of restricted three body problem which was first time studied by \cite{Chermnykh1987VeLen}. \cite{PapadakisKanavos2007Ap&SS}  given numerical exploration of Chermnykh's problem, in which the equilibrium points and zero velocity curves studied numerically also the non-linear stability for the triangular Lagrangian points are computed numerically for the Earth-Moon and Sun-Jupiter mass distribution when the angular velocity varies. The mass reduction factor  $q_1=1-\frac{F_p}{F_g}$   is expressed in terms of the particle radius $\mathbf{a}$, density $\rho$ radiation pressure efficiency factor $\chi$ (in C.G.S. system):\( q_1=1-\frac{5.6\times{10^{-5}}}{\mathbf{a}\rho}\chi
\). Where $F_p$ is solar radiation pressure force which is  exactly apposite to the gravitational attraction force $F_g$ and change with the distance by the same law it is possible to consider that the result of action of this force will lead to reducing the effective mass of the  Sun or particle.

~\cite{KushvahBR2006} examined the linear stability of triangular
equilibrium points in the generalized photogravitational restricted
three body problem with Poynting-Robertson drag, $L_4$ and $L_5$
points became unstable due to P-R drag which is very remarkable and
important, where as they are linearly stable in classical problem when $0<\mu<\mu_{Routh}=0.0385201$. Further the normalizations  of Hamiltonian and nonlinear stability of $L_{4(5)}$ in the present of P-R drag has been studied by \cite{Kushvah2007BASI,Kushvah2007Ap&SS.312,Kushvah2007EM&P..101} 

Normal forms are a standard tool in Hamiltonian mechanics
to study the dynamics in a neighbourhood of invariant objects. Usually,
these normal forms are obtained as divergent series, but their
asymptotic character is what makes them useful. From theoretical
point of view, they provide nonlinear approximations to the dynamics
in a neighbourhood of the invariant object, that allows to obtain
information about the real  solutions of the system by taking the
normal form upto a suitable finite order. In this case, it is well
known that under certain(generic)  non-resonance conditions, the
remainder of this finite normal form turns out to be exponentially
small with respect to some parameters. Those series are usually
divergent  on open sets, it is still possible in some cases to prove
convergence on certain sets with empty interior(Cantor-like sets) by
replacing the standard linear normal form scheme by a quadratic
one. 

From a more practical point of view, normal forms can be used as a
computational method to obtain very accurate approximations to the
dynamics in a neighbourhood of the selected invariant object, by
neglecting the remainder. They have been applied, for example, to
compute invariant manifolds or invariant  tori. To do that, it is
necessary to compute the explicit expression of the normal form and
of the (canonical) transformation that put the Hamiltonian into this
reduced form. A context where this computational formulation has
special interest is in some celestial mechanics models, that can be
used to approximate the dynamics of some real world problems.

The linearization of Hamiltonian, have been  studied numerically and analytically. We have examined the effect of  gravitational potential from the belt, oblateness effect and radiation effect.

\section{Equations of Motion and Position of Equilibrium Points}
\begin{center}
\begin{figure}
    \includegraphics[scale=.6]{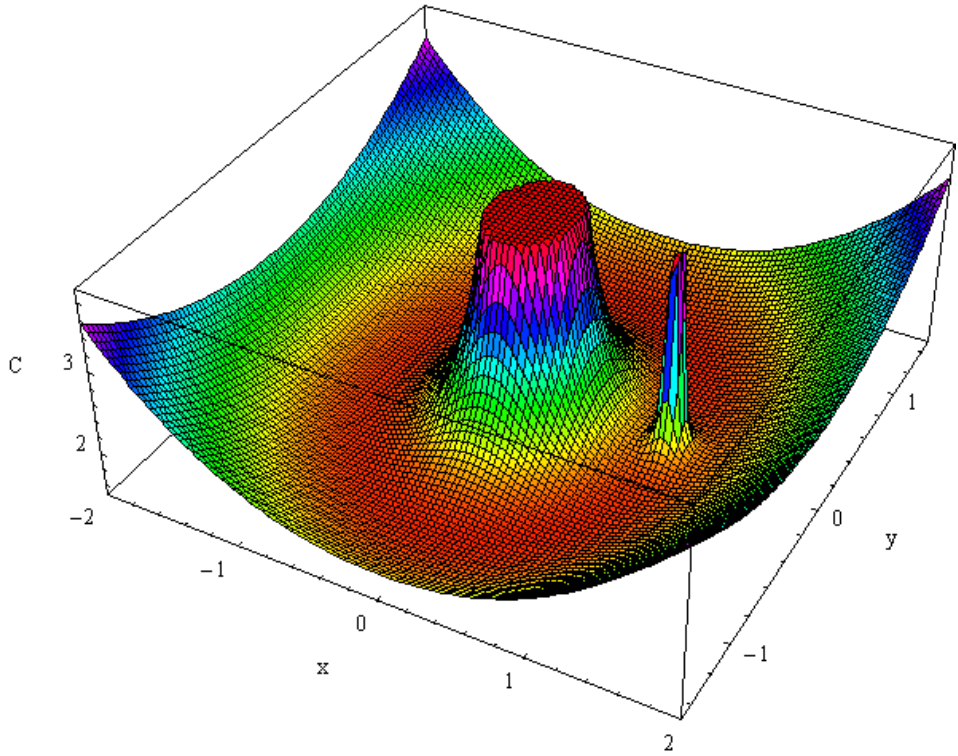}  \includegraphics[scale=.85]{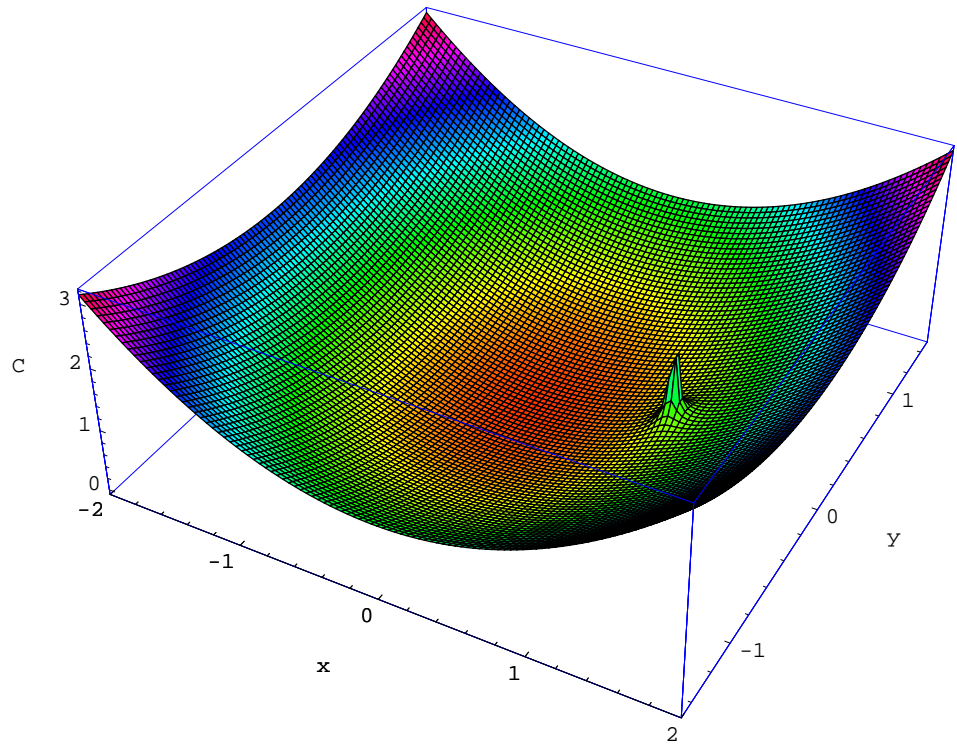} 
\caption{The  zero velocity surfaces for $\mu=0.025,r_c=0.9999,T=0.01$, for all values of $A_2, M_b $}
 \label{fig:plotzvc}
\end{figure}
\end{center}

\begin{center}
  \begin{figure}
    \includegraphics[scale=.8]{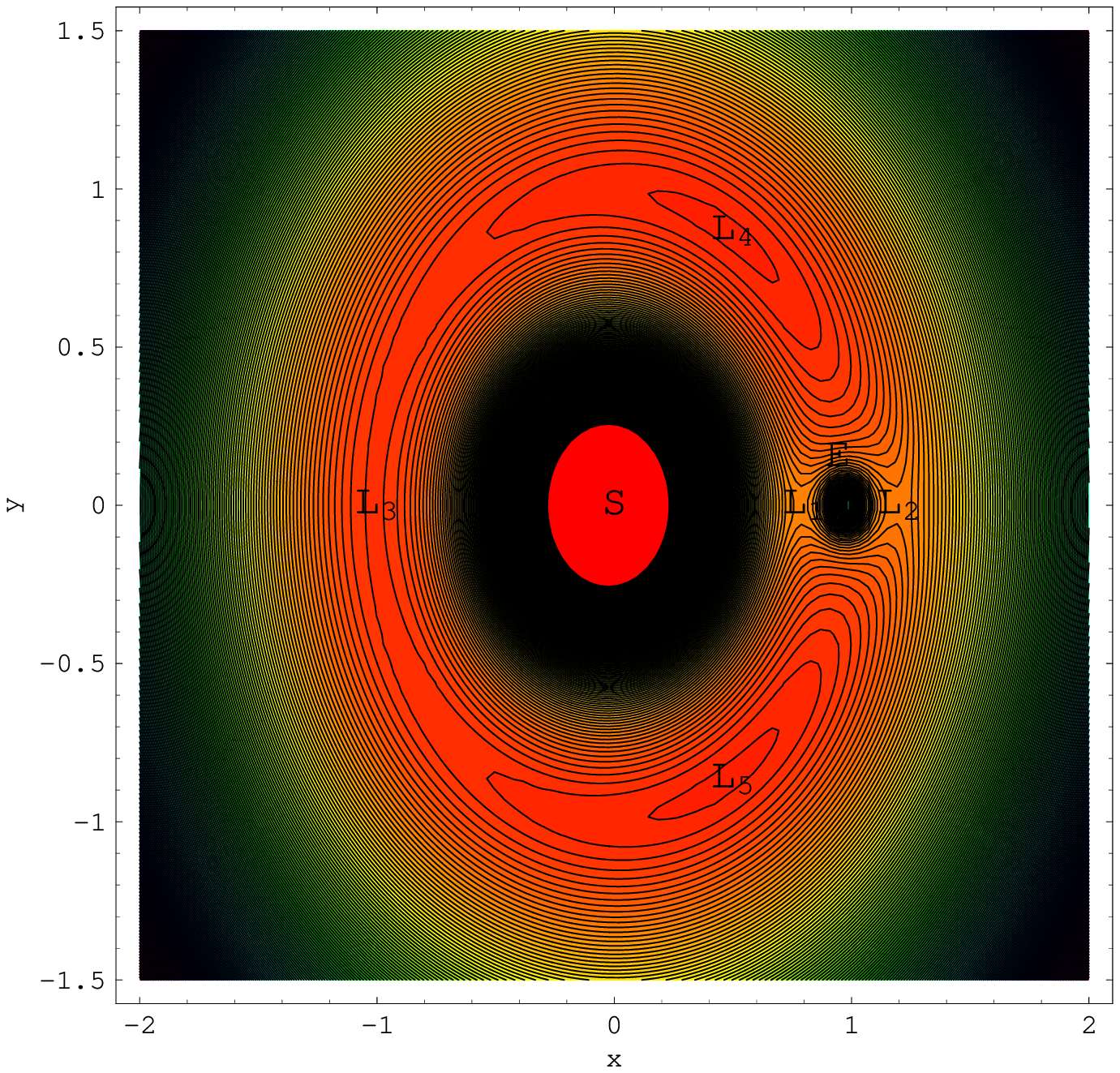}\\
 \includegraphics[scale=.75]{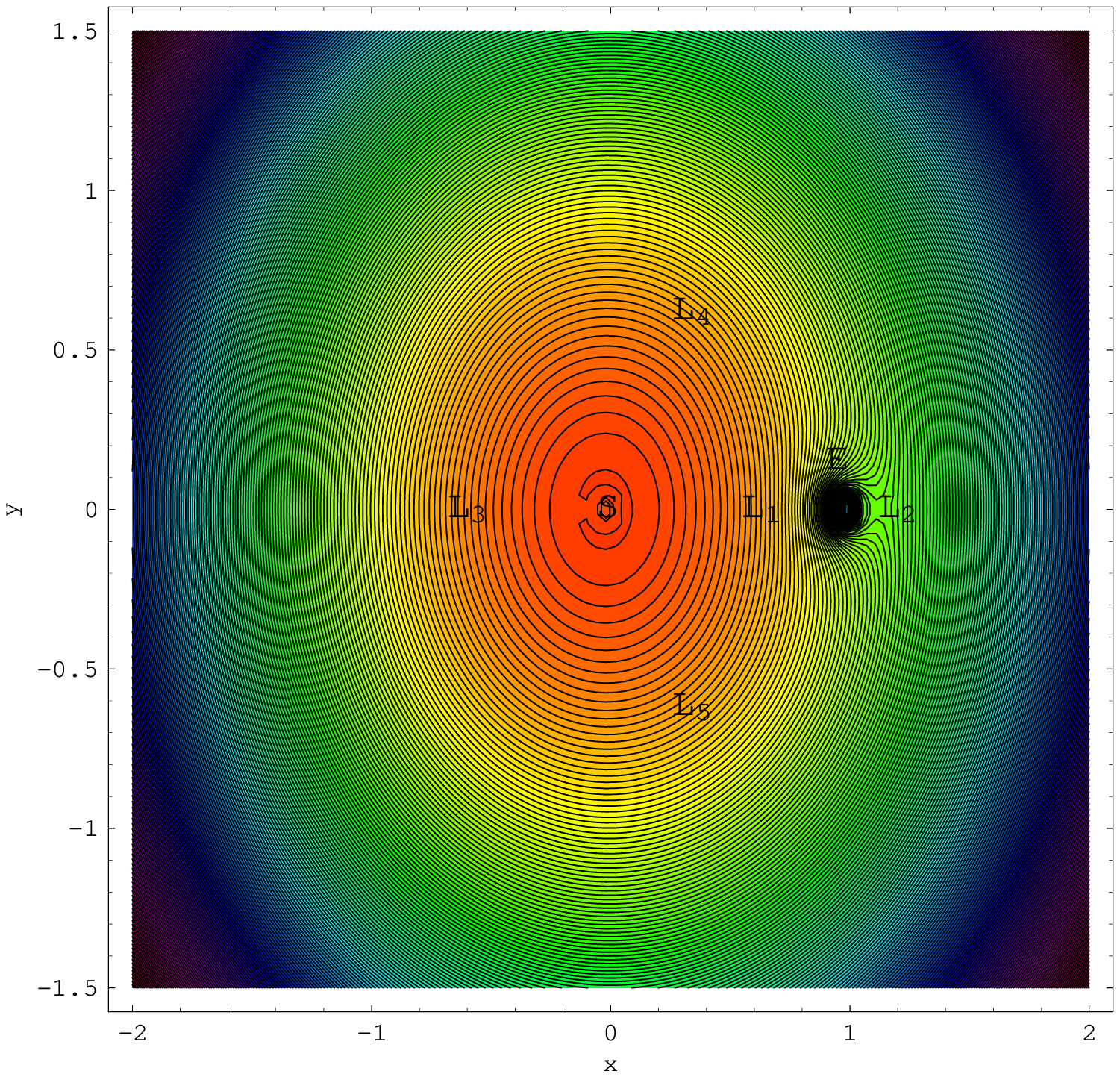} 
   \caption{The figure show the zero velocity curves for $\mu=0.025, T=0.01$, in frame one $q_1=0.00$, in  second frame  show the curves in classical case.}
  \label{fig:jvc&q0}
\end{figure}
\end{center}

Let us consider the model proposed by \cite{MiyamotoNagai1975PASJ}, according to this model the potential of belt is given by:
\begin{eqnarray}
 V (r,z)=\frac{\mathbf{b}^2M_b\left[\mathbf{a}r^2 +\left(\mathbf{a}+3N\right)\right]\left( \mathbf{a}+N\right)^2}{N^3\left[r^2+\left(\mathbf{a}+N\right)^2\right]^{5/2}}
\end{eqnarray}
where $M_b$ is the total mass of the belt and $r^2=x^2+y^2$, $\mathbf{a,b}$ are parameters which determine the density profile of the belt, if  $\mathbf{a}=\mathbf{b}=0$ then the  potential equals to the one by a point mass. The parameter $\mathbf{a}$  \lq\lq flatness parameter\rq\rq and $\mathbf{b}$  \lq\lq core parameter\rq\rq, 
 where $N=\sqrt{z^2+\mathbf{b}^2}$,\ $T=\mathbf{a}+\mathbf{b}$, $z=0$.
Then we obtained \begin{equation}
 V(r,0)=-\frac{M_b}{\sqrt{r^2+T^2}}\label{eq:Vr0}
\end{equation}
and  \(V_x=\frac{M_b x}{\left(r^2+T^2\right)^{3/2}},\ V_y=\frac{M_b y}{\left(r^2+T^2\right)^{3/2}} \).
As in   \cite{Kushvah2008Ap&SS,Kushvah2008Ap&SS.315}, we consider the barycentric rotating co-ordinate system $Oxyz$ relative to inertial system with angular velocity $\omega$ and common $z$--axis.  We have taken line joining the primaries as $x$--axis. Let $m_1, m_2$ be the masses of bigger primary(Sun)  and smaller primary(Earth) respectively. Let  $Ox$, $Oy$  in the equatorial plane of smaller primary  and $Oz$ coinciding with the polar axis of $m_{2}$. Let $r_{e}$, $r_{p}$ be the equatorial and polar radii of $m_{2}$ respectively,  $r$ be the distance between primaries.  Let infinitesimal mass $m$ be placed at the  point $P(x,y,0)$. We take units such that sum of the masses and distance between primaries as  unity, the unit of time taken such that the Gaussian constant of gravitational $\Bbbk^{2}=1$. Then perturbed mean motion $n$ of the primaries is given by  $n^{2}=1+\frac{3A_{2}}{2}+\frac{2M_b r_c}{\left(r_c^2+T^2\right)^{3/2}}$, where $r_c^2=(1-\mu)q_1^{2/3}+\mu^2$, $A_{2}=\frac{r^{2}_{e}-r^{2}_{p}}{5r^{2}}$ is oblateness coefficient of $m_{2}$.
where $\mu=\frac{m_{2}}{m_{1}+m_{2}}$  is mass parameter,  $1-\mu=\frac{m_{1}}{m_{1}+m_{2}}$ with $m_{1}>m_{2}$. Then coordinates of $m_1$ and $m_2$ are  $(x_1,0)=(-\mu,0)$ and  $(x_2,0)=(1-\mu,0)$ respectively. The mass of Sun $m_1\approx 1.989\times 10^{30}kg$ $\approx 332,946m_2$(The mass of Earth), hence mass parameter for this system is $\mu=3.00348\times10^{-6}$.  In the above mentioned reference system, we determine  the equations of motion of the infinitesimal mass particle in $x y$-plane.  Now using  \cite{MiyamotoNagai1975PASJ} profile and  \cite{Kushvah2008Ap&SS,Kushvah2008Ap&SS.315}, the equations of motion are given  by:
\begin{eqnarray}
\ddot{x}-2n\dot{y}&=&U_{x}-V_x=\Omega_x ,\label{eq:Omegax}\\
\ddot{y}+2n\dot{x}&=&U_{y}-V_y=\Omega_y\label{eq:Omegay}
 \end{eqnarray}
where
\begin{eqnarray}
&&\Omega_x=n^{2}x-\frac{(1-\mu)q_1(x+\mu)}{r^3_1}-\frac{\mu(x+\mu-1)}{r^3_2}-\frac{3}{2}\frac{\mu{A_2}(x+\mu-1)}{r^5_2}-\frac{M_b x}{\left(r^2+T^2\right)^{3/2}},\nonumber\\
&&\Omega_y=n^{2}y
-\frac{(1-\mu)q_{1}{y}}{r^3_1}
-\frac{\mu{y}}{r^3_2}-\frac{3}{2}\frac{\mu{A_2}y}{r^5_2}-\frac{M_b y}{\left(r^2+T^2\right)^{3/2}},\nonumber\\
&&\Omega=\frac{n^2(x^2+y^2)}{2}+\frac{(1-\mu)q_1}{r_1}+\frac{\mu}{r_2}+\frac{\mu
 A_2}{2r_2^3}+\frac{M_b}{\left(r^2+T^2\right)^{1/2}}\label{eq:OmegaFF}
 \end{eqnarray}
The energy integral of the problem is given by  $C=2\Omega-{\dot{x}}^2-{\dot{y}}^2$, where the quantity $C$ is the Jacobi's  constant. The zero velocity curves $C=2\Omega(x,y)\label{eq:C}$ are presented in  figure (~\ref{fig:plotzvc})  for the entire range of parameters$ A_2, M_b$ and $ q_1=1,0$.  We have seen that  there are closed curves around the $L_{4(5)}$ so they are stable but the stability range reduced(~\ref{fig:jvc&q0})] due radiation  effect.  The closed curves around $L_{4(5)}$ disappeared when $q_1=0$.

The  position of equilibrium points  are given by putting  $\Omega_x=\Omega_y=0$ i.e.,  
\begin{eqnarray}
&& n^{2}x-\frac{(1-\mu)q_1(x+\mu)}{r^3_1}-\frac{\mu(x+\mu-1)}{r^3_2}\nonumber\\&&-\frac{3}{2}\frac{\mu{A_2}(x+\mu-1)}{r^5_2}-\frac{M_b x}{\left(r^2+T^2\right)^{3/2}} =0\label{eq:eq1pts},\\
&&n^{2}y
-\frac{(1-\mu)q_{1}{y}}{r^3_1}-\frac{\mu{y}}{r^3_2}-\frac{3}{2}\frac{\mu{A_2}y}{r^5_2}\nonumber\\&&-\frac{M_b y}{\left(r^2+T^2\right)^{3/2}} =0\label{eq:eq2pts}\end{eqnarray}
 from equations  (~\ref{eq:eq1pts}, ~\ref{eq:eq2pts}) we obtained:
\begin{eqnarray}
r_1&&=q_1^{1/3}\left[1-\frac{A_2}{2}+\frac{(1-2r_c)M_b \left(1-\frac{3\mu A_2}{2(1-\mu)}\right)}{3\left(r_c^2+T^2\right)^{3/2}}\right]\label{eq:r1},\\ 
r_2&&=1+\frac{\mu(1-2r_c)M_b}{3\left(r_c^2+T^2\right)^{3/2}} \label{eq:r2}
\end{eqnarray}
From above, we obtained:
\begin{eqnarray} &&x=-\mu\pm\left[\left(\frac{q_1}{n^2}\right)^{2/3}\left[1+\frac{3A_2}{2}\right.\right.\nonumber\\&&\left.\left.-\frac{(1-2r_c)M_b \left(1-\frac{3\mu A_2}{2(1-\mu)}\right)}{\left(r_c^2+T^2\right)^{3/2}}\right]^{-2/3}-y^2\right]^{1/2} \label{eq:x1lag}\\
&&x=1-\mu\nonumber\\&&\pm\Bigl[\left[1-\frac{\mu(1-2r_0)M_b}{\left(r_c^2+T^2\right)^{3/2}}\right]^{-2/3}-y^2\Bigr]^{1/2}\label{eq:x2lag}
\end{eqnarray}
 
The triangular equilibrium points are given by putting  $\Omega_x=\Omega_y=0$, $y\neq{0}$,   then from equations (~\ref{eq:Omegax}) and (~\ref{eq:Omegay}) we obtained the triangular equilibrium points as:%
\begin{eqnarray}
&&x=-\mu+\frac{q_1^{2/3}}{2}(1-A_2)+ \frac{(1-2r_c) M_b\left[\left\{1-\frac{3\mu A_2}{(1-\mu)}\right\}q_1^{2/3}-1\right]}{3\left(r_c^2+T^2\right)^{3/2}}\label{eq:xl4}\end{eqnarray}
\begin{eqnarray}
 &&y=\pm \frac{q_1^{2/3}}{2} \left[4-q_1^{2/3}+2\left(q_1^{2/3}-2\right)A_2\right.\nonumber\\&&\left.-\frac{4(2r_c-1)M_b\left[\left\{\left(q_1^{2/3}-3\right)-\frac{3\mu A_2\left(q_1^{2/3}-3\right)}{2(1-\mu)}\right\}\right]}{3\left(r_c^2+T^2\right)^{3/2}}\right]^{1/2}\label{eq:yl4} 
 \end{eqnarray}
All these results are similar with \cite{Szebehely1967}, \cite{Ragosetal1995},   \cite{Kushvah2008Ap&SS,Kushvah2008Ap&SS.315} and others. 
\section{Linearization of the Hamiltonian}
\label{sec:1storder}
We have to expand the Lagrangian function in power series of $x$ and $y$, where $(x, y)$ are the coordinates of the triangular equilibrium points. We will  examine the stability of the triangular equilibrium points. For this we will utilize the method of \cite{Whittaker1965}. By taking $H_2$, we will  consider  linear equations this we have established the relations between perturbed basic frequencies. The Lagrangian function of the problem can be written as
\begin{eqnarray}
L&=&\frac{1}{2}(\dot{x}^2+\dot{y}^2)+n(x\dot{y}-\dot{x}y)+\frac{n^2}{2}(x^2+y^2)\nonumber\\&&+\frac{\amc{q_1}}{r_1}+\frac{\mu}{r_2}+\frac{\mu{A_2}}{2r^3_2}+\frac{M_b}{\left(r^2+T^2\right)^{1/2}}\label{eq:L}
\end{eqnarray}
and the Hamiltonian $H=-L+p_x\dot{x}+p_y\dot{y}$, where $p_x,p_y$
are the momenta coordinates given by \[
p_x=\frac{\partial{L}}{\partial{\dot{x}}}=\dot{x}-ny,
\quad
p_y=\frac{\partial{L}}{\partial{\dot{y}}}=\dot{y}+nx
\]
Let us  suppose  $q_1=1-\epsilon$, with $|\epsilon|<<1$ and rewriting the  coordinates of triangular equilibrium points as: 
\begin{eqnarray}
x&=\frac{\gamma}{2}-\frac{\epsilon}{3}-\frac{A_2}{2}+\frac{A_2
\epsilon}{3}+\frac{2M_b \epsilon}{9}-\frac{\mu A_2 M_b}{2(1-\mu)}\left(1-\frac{2}{3}\epsilon\right) \\
y&=\pm \frac{\sqrt{3}}{2}\Bigl\{1-\frac{5\epsilon}{9}-\frac{A_2}{3}-\frac{2A_2
\epsilon}{9}-\frac{4M_b}{9}-\frac{8M_b
\epsilon}{27}+\frac{\mu A_2 M_b\epsilon}{9(1-\mu)}\Bigr\}
\end{eqnarray}
where $\gamma=1-2\mu$.
 We shift the origin to $L_4$, change
$x\rightarrow {x_*}+x$,  $y\rightarrow{y_*}+y$. and  $a=x_*+\mu,
b=y_*$,then 
\begin{eqnarray}
a&= \frac{1}{2} \left\{1-\frac{2\epsilon}{3}-A_2+\frac{2A_2
\epsilon}{3}+\frac{4M_b \epsilon}{9}-\frac{\mu  A_2 M_b}{(1-\mu)}\left(1-\frac{2}{3}\epsilon\right) \right\}\label{eq:a}\\
b&=\frac{\sqrt{3}}{2}\left\{1-\frac{2\epsilon}{9}-\frac{A_2}{3}-\frac{2A_2
\epsilon}{9}-\frac{4M_b}{9}-\frac{8M_b
\epsilon}{27}+\frac{\mu A_2 M_b\epsilon}{9(1-\mu)}\right\}\label{eq:b}
\end{eqnarray}
We have to expand $L$ in power series of $x $ and $y$, for this we
use Taylor's expansion i.e.
\begin{eqnarray}
&&f(x,y)=f(0,0,)+[x f_x (0,0)+y f_y
(0,0)]+\frac{1}{2}[x^2f_{xx}(0,0)+2xyf_{xy}(0,0)+y^2f_{yy}(0,0)]
\nonumber\\&&+\frac{1}{3!}[x^3f_{xxx}(0,0)+3
x^2yf_{xxy}(0,0)+3xy^2f_{xyy}(0,0)+y^3f_{yyy}(0,0)]+\dots \label{eq:Texp}
\end{eqnarray}
Using  these values in (~\ref{eq:Texp}) and with the help of   (~\ref{eq:a}) and  (~\ref{eq:b})  we get
\begin{eqnarray}
 L&=&L_0+L_1+L_2+L_3+\cdots \\
H&=&H_0+H_1+H_2+H_3+\cdots =-L+p_x{\dot{x}}+p_y{\dot{y}}
  \end{eqnarray}
  where $L_0,L_1,L_2,\ldots$ are
\begin{eqnarray}
L_0&=&\left[(a-\mu)^2 +b^2\right]n^2+\frac{(1-\mu) q_1 }{\sqrt{a^2+b^2}}\nonumber\\&&+\frac{\mu  }{\sqrt{(a-1)^2+b^2}}\left\{1+\frac{A_2}{2(a-1)^2+b^2}\right\}+\frac{M_b}{\left((a-\mu)^2+b^2+T^2\right)^{3/2}}
\end{eqnarray}
\begin{eqnarray}
L_1&=&x\left[2 n^2 (a-\mu)-\frac{a (1-\mu)q_1}{\left(a^2+b^2\right)^{3/2}}-\frac{3(a-\mu)M_b}{\left((a-\mu)^2+b^2+T^2\right)^{5/2}}\right.\nonumber\\&&\left.-\frac{(a-1)\mu  }{\left(a-1)^2+b^2\right)^{3/2}}\left\{1+\frac{3 A_2}{2(a-1)^2+b^2}\right\}\right]\nonumber\\&& y\left[2b n^2-\frac{b  (1-\mu)q_1}{\left(a^2+b^2\right)^{3/2}}-\frac{3bM_b}{\left((a-\mu)^2+b^2+T^2\right)^{5/2}}\right.\nonumber\\&&\left.-\frac{b\mu  }{\left(a-1)^2+b^2\right)^{3/2}}\left\{1+\frac{3 A_2}{2(a-1)^2+b^2}\right\}\right]
\end{eqnarray}
\begin{eqnarray}
L_2&=&\frac{(\dot x^2+ \dot y^2)}{2}+n(x\dot y-\dot x y)+
\frac{n^2}{2}(x^2+y^2)-Ex^2-Fy^2-G xy
\end{eqnarray}
\begin{eqnarray}
E&=&\frac{1}{1728}\left[ 2-324(2-25\mu)M_b-\frac{(2r_c-1)\left\{2+15(2-7\mu)M_b\right\}M_b}{\left(r_c^2+T^2\right)^{3/2}}\right.\nonumber\\&+&2\epsilon\left\{-54-270\mu +135(10-81\mu)M_b\right.\nonumber\\&-&\left.\frac{(2r_c-1)\left\{146-33\mu+15(118-205\mu)M_b\right\}M_b}{\left(r_c^2+T^2\right)^{3/2}}\right\}\nonumber\\&+&6A_2\left\{162-432\mu +135(2+39\mu)M_b\right.\nonumber\\&+&\left.\frac{(2r_c-1)\left\{146-240\mu-15(10-259\mu)M_b\right\}M_b}{\left(r_c^2+T^2\right)^{3/2}}\right\}\nonumber\\&+&\epsilon A_2\left\{144-5022\mu -270(4-395\mu)M_b\right.\nonumber\\&&-\left.\left.\frac{(2r_c-1)\left\{458+540\mu+15(1926-7399\mu)M_b\right\}M_b}{\left(r_c^2+T^2\right)^{3/2}}\right\}\right]\label{eq:E}
\end{eqnarray}
  \begin{eqnarray}
F&=&\frac{-1}{576}\left[360+108(22+85\mu)M_b+\frac{(2r_c-1)\left\{2+5(6+35\mu)M_b\right\}M_b}{\left(r_c^2+T^2\right)^{3/2}}\right.\nonumber\\&+&2\epsilon\left\{54-18\mu +45(62+309\mu)M_b\right.\nonumber\\&-&\left.\frac{(2r_c-1)\left\{10+39\mu+5(202+1565\mu)M_b\right\}M_b}{\left(r_c^2+T^2\right)^{3/2}}\right\}\nonumber\\&+&6A_2\left\{126+45(26+99\mu)M_b\right.\nonumber\\&+&\left.\frac{(2r_c-1)\left\{46+5(142+791\mu)M_b\right\}M_b}{\left(r_c^2+T^2\right)^{3/2}}\right\}\nonumber\\&+&\epsilon A_2\left\{576-18\mu+90(228+1147\mu)M_b\right.\nonumber\\&&+\left.\left.\frac{(2r_c-1)\left\{522+526\mu+5(3834+27923\mu)M_b\right\}M_b}{\left(r_c^2+T^2\right)^{3/2}}\right\}\right]\label{eq:F}
  \end{eqnarray}
   \begin{eqnarray}
G&=&\frac{-1}{288\sqrt{3}}\left[ 648-1296\mu+1620(2+3\mu)M_b\right.\nonumber\\&+&\left.\frac{12(2r_c-1)\left\{22-44\mu+5(34+39\mu)M_b\right\}M_b}{\left(r_c^2+T^2\right)^{3/2}}\right.\nonumber\\&+&2\epsilon\left\{198-666\mu +45(94+201\mu)M_b\right.\nonumber\\&+&\left.\frac{2(2r_c-1)\left\{190-423\mu+135(18+65\mu)M_b\right\}M_b}{\left(r_c^2+T^2\right)^{3/2}}\right\}\nonumber\\&+&6A_2\left\{126-468\mu +45(26+15\mu)M_b\right.\nonumber\\&+&\left.\frac{(2r_c-1)\left\{96-260\mu+15(66+133\mu)M_b\right\}M_b}{\left(r_c^2+T^2\right)^{3/2}}\right\}\nonumber\\&+&\epsilon A_2\left\{1368-4806\mu +90(280+429\mu)M_b\right.\nonumber\\&+&\left.\left.\frac{(2r_c-1)\left\{2106-655\mu+35(982+3025\mu)M_b\right\}M_b}{\left(r_c^2+T^2\right)^{3/2}}\right\}\right]\label{eq:G}  \end{eqnarray}
\subsection{Perturbed Basic Frequencies}
\label{subsec:pertbasic}
Using the method in  \cite{Whittaker1965}, to find the canonical transformation from the phase space $(x,y,p_x,p_y)$ into the phase space product of the angle co-ordinates $(\phi_1,\phi_2)$ and the action momenta co-ordinates $I_1,I_2$ and of the first order in $I_1^{1/2},I_2^{1/2}$. We consider the following  linear equations in the variables $x, y$:
 \begin{equation}
 \begin{array}{l c l}
 -\lambda p_x& = & \frac{\partial{H_2}}{\partial x},\\&&\\
 -\lambda p_y& = & \frac{\partial{H_2}}{\partial y},\\
\quad AX&=&0
 \end{array} \quad
 \begin{array}{l c l }
 \lambda x& = & \frac{\partial{H_2}}{\partial p_x},\\&&\\
 \lambda y& = & \frac{\partial{H_2}}{\partial p_y},\\&&
 \label{eq:AX}\end{array}
  \end{equation}
  \begin{equation}
  X=\left[\begin{array}{c}
  x\\
  y\\
  p_x\\
  p_y \end{array}\right] \quad {and}
 \quad
 A=\left[\begin{array}{c c c c}
  2E & G&\lambda& -n\\
G&2F&n&\lambda\\
  -\lambda& n& 1& 0\\
  -n & -\lambda& 0& 1\end{array}\right]
  \end{equation}
The characteristic equation  of the  Hamiltonian
 \[
H_2=\frac{p_x^2+p_y^2}{2}+n(yp_x-xp_y)+Ex^2+Fy^2+Gxy
\] is given by $|A|=0$ i.e.
 \begin{equation}
 \lambda^4+2(E+F+n^2)\lambda^2+4EF -G^2+n^4-2n^2(E+F)=0 \label{eq:ch}
 \end{equation}
 Stability is assured  only when the discriminant  $D>0$,  where
  \begin{equation}
 D=4(E+F+n^2)^2-4\bigl\{4EF-G^2+n^4-2n^2(E+F)\bigr\}
 \end{equation}
 from stability condition we obtained:
  \begin{eqnarray}
 \mu&<&\mu_{c_0}+ 0.125885\epsilon+M_b\left\{0.571136+1.73097\epsilon\right.\nonumber\\&&+ \frac{ 0.219964-0.398363 \epsilon+r_c(0.202129+0.114305\epsilon)}{\left(r_c^2+T^2\right)^{3/2}}\left.\right\}\nonumber\\&&-\left[0.0627796 -0.112691\epsilon+M_bA_2\left\{0.281354-1.53665\epsilon\right.\right.\nonumber\\&&\left.\left.+ \frac{0.654936-0.669428\epsilon+r_c(0.195486+0.350878\epsilon)}{\left(r_c^2+T^2\right)^{3/2}}\right\}\right]\end{eqnarray}
 where $\mu_{c_0}=0.0385209$.
 When $D>0$ the roots $\pm i \omega_1$ and $\pm i \omega_2$ ($\omega_1,\omega_2$ being the long/short -periodic frequencies) are related to each other  as
 \begin{eqnarray}
  & &\omega_1^2+\omega_2^2= \frac{1}{54}\left[27((1+\mu)\epsilon-2)+9[-18+36\mu+(22+69\mu)]\right.\nonumber\\&&\left.+\frac{81M_b}{2}(12+30\epsilon+30\mu+95\mu\epsilon)+135M_bA_2(18+58\epsilon+45\mu+188\mu\epsilon)\right.\nonumber\\&&\left.+\frac{1}{2\left(r_c^2+T^2\right)^{3/2}}\Bigl\{M_b\left(180+2(2r_c-1)(44+21)\mu\epsilon\right.\right.\nonumber\\&&+\left.\left.72r_c+4r_cA_2[78+180\mu+(253+873\mu)\epsilon]\right)\Bigr\}\right]\label{eq:w1+w2}\end{eqnarray}
\begin{eqnarray}
   &&\omega_1^2\omega_2^2= -\frac{1}{72}\left[108[\epsilon+3M_b(2+7\epsilon)]+18[31\epsilon+M_b(144+647\epsilon)]A_2\right.\nonumber\\&&+\left.\frac{M_b\left\{4 [36-47\epsilon+4r_c(9+28\epsilon)]+3[74-297\epsilon+4r_c(44+175\epsilon)]A_2\right\}}{\left(r_c^2+T^2\right)^{3/2}}\right]\nonumber\\&&+\mu\left[\frac{27}{4}+\frac{99\epsilon}{8}+\frac{117}{4}+73A_2\epsilon +M_b\left\{\frac{45}{4}+81A_2+\frac{273\epsilon}{8}+\frac{2357A_2\epsilon }{8}\right.\right.\nonumber\\&&\left.\left.+\frac{396(2r_c-1)+4(772r_c-395)\epsilon+12(326r_c-181)A_2+(18452r_c-9667)A_2\epsilon}{72\left(r_c^2+T^2\right)^{3/2}} \right\}\right]\nonumber\\&&+\mu^2\left[\frac{-27}{4}+\frac{111\epsilon}{8}+\frac{117}{4}+\frac{161A_2\epsilon}{2} +M_b\left\{\frac{405}{4}+\frac{495A_2}{2}+\frac{4185\epsilon}{16}+\frac{25275A_2\epsilon }{16}\right.\right.\nonumber\\&&\left.\left.+\frac{(2r_c-1)[198-838\epsilon-1014A_2-5117A_2\epsilon]}{36\left(r_c^2+T^2\right)^{3/2}} \right\}\right]\label{eq:w1w2}\end{eqnarray}
where  $\omega_j (j=1,2)$ satisfy the property $(0<\omega_2<\frac{1}{\sqrt{2}}<\omega_1<1)$.
\begin{figure}
\includegraphics[scale=.9]{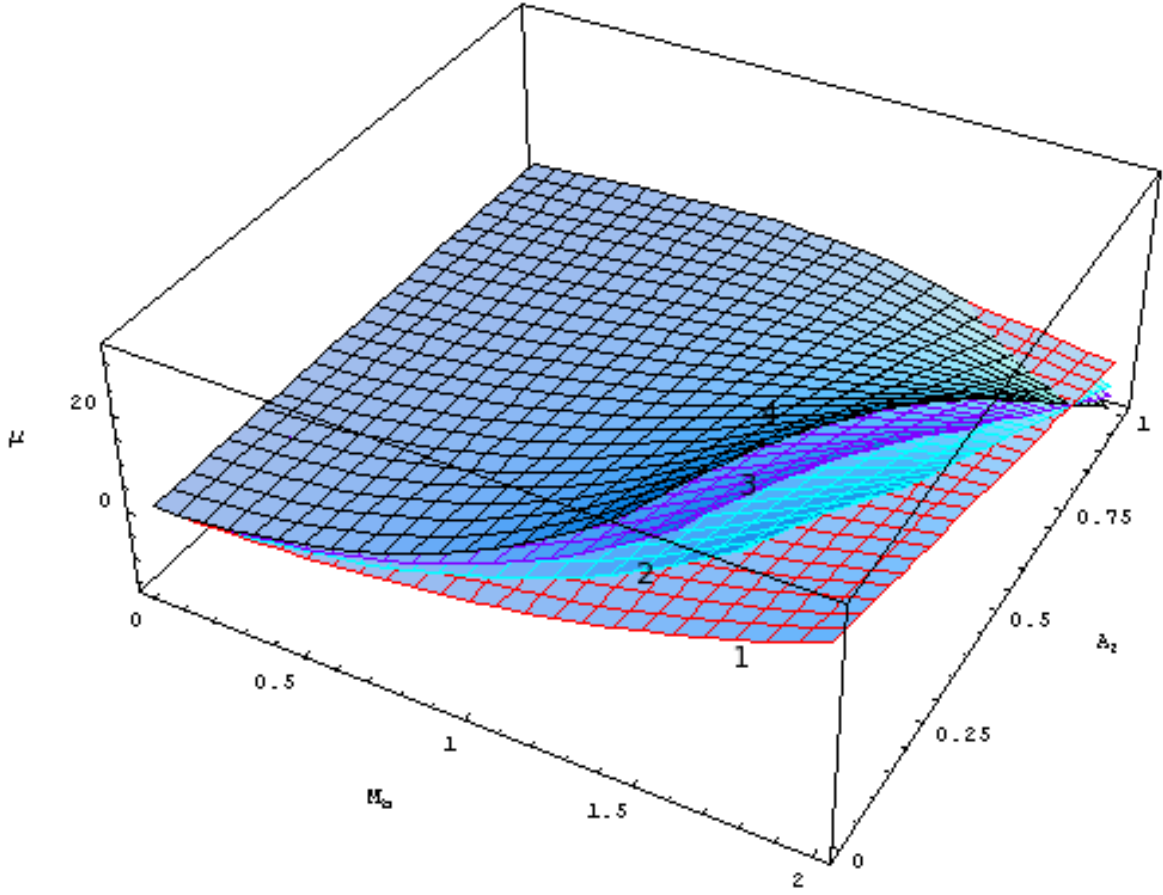} 
   \caption{The  critical mass ratios  in the 3D plot for $M_b-A_2-\mu$ and $q_1=1, 0.75, 0.5, 0.25$, $r_c=0.9999, T=0.01$}
  \label{fig:mu}
\end{figure}
The mass parameter is shown in figure $0\leq M_b\leq2$, $0\leq A_2\leq 1$, surface (1)$q_1=1$ (2) $q_1=0.75$, (3)$q_1=0.5$,(4)$q_1=0.25$.
\subsection{The Normal Coordinates}
\label{subsec:detnorm} For expressing $H_2$ in a simpler form, we
consider the set of linear equations(~\ref{eq:AX}) the solution of
which can be obtained as
\begin{eqnarray}
\frac{x}{(2n\lambda- G)}&=&\frac{y}{(\lambda^2-n^2+2E)}=\frac{p_x}{(n\lambda^2-G\lambda-2nE+n^3)}\nonumber\\&&=\frac{p_y}{(\lambda^3+n^2\lambda+2E\lambda-nG)} \end{eqnarray}
Substituting $\lambda=\pm i\omega_1$ and $\pm i\omega_2$, we obtain the solution sets as
\[x_j=K_j(2ni\omega_j-G),\quad p_{x,j}=K_j(-n\omega_j^2-iG\omega_j-2En+n^3)\]
\[y_j=K_j(-\omega_j-n^2+2E), \quad p_{y,j}=K_j\left\{-i\omega_j^3+i\omega_j(n^2+2E)-Gn\right\}\]
\[x_{j+2}=K_{j+2}(2ni\omega_j-G), \quad p_{x,j+2}=K_{j+2}(-n\omega_j^2+iG\omega_j-2En+n^3)\]
\[y_{j+2}=K_{j+2}(-\omega_j^2-n^2+2E), \quad p_{y,j+2}=K_{j+2}\left\{\omega_j^3-i\omega_j(n^2+2E)-Gn\right\}\]
where $j=1,2$ and $K_j,K_{j+2}$ are constants of proportionality. Following the method for reducing $H_2$ to the normal form, as in \cite{Whittaker1965}, use the transformation
\begin{equation}
 X=JT \quad  {where}  \quad X=\left[\begin{array}{c}
x\\y\\p_x\\p_y\end{array}\right],\quad T=\left[\begin{array}{c}
Q_1\\Q_2\\P_1\\P_2\end{array}\right]\label{eq:XJT}
\end{equation}
and $J=[J_{ij}]_{1\leq i, j \leq 4}$
i.e.
\begin{equation}
  J=\left[\begin{array}{c c c c}
  x_1-i\frac{\omega_1 x_3}{2}& -x_2-i\frac{\omega_2 x_4}{2}&-i\frac{x_1}{\omega_1}+\frac{x_3}{2}&-i\frac{x_2}{\omega_2}+\frac{x_4}{2}\\
  y_1-i\frac{\omega_1 y_3}{2}& -y_2-i\frac{\omega_2 y_4}{2}&-i\frac{y_1}{\omega_1}+\frac{y_3}{2}&-i\frac{y_2}{\omega_2}+\frac{y_4}{2}\\
 p_{x,1}-i\frac{\omega_1 p_{x,3}}{2}& -p_{x,2}-i\frac{\omega_2 p_{x,4}}{2}&-i\frac{p_{x,1}}{\omega_1}+\frac{p_{x,3}}{2}&-i\frac{p_{x,2}}{\omega_2}+\frac{p_{x,4}}{2}\\
  p_{y,1}-i\frac{\omega_1 p_{y,3}}{2}& -p_{y,2}-i\frac{\omega_2 p_{y,4}}{2}&-i\frac{p_{y,1}}{\omega_1}+\frac{p_{y,3}}{2}&-i\frac{p_{y,2}}{\omega_2}+\frac{p_{y,4}}{2} \end{array}\right] \end{equation}
 $P_i= (2 I_i\omega_i)^{1/2}\cos{\phi_i},  \quad Q_i= (\frac{2
I_i}{\omega_i})^{1/2}\sin{\phi_i}, \quad (i=1,2)$
Under normality conditions:
\begin{eqnarray*}
x_1p_{x,3}-x_3p_{x,1}+y_1p_{y,3}-y_3p_{y,1}&=&1\\
x_2p_{x,4}-x_4p_{x,4}+y_2p_{y,4}-y_4p_{y,2}&=&1
\end{eqnarray*}
Equivalently,
\begin{eqnarray}
&&-4i\omega_1K_1K_2\bigl\{\omega_1^2(F-E+n^2)\nonumber\\&&+G^2+2E^2+3n^2E+n^2F-2EF-2n^4\bigr\}=1\label{eq:4ik1}\\&&
-4i\omega_2K_2K_4\bigl\{\omega_2^2(F-E+n^2)\nonumber\\&&+G^2+2E^2+3n^2E+n^2F-2EF-2n^4\bigr\}=1\label{eq:4ik2}
\end{eqnarray}
$K_j's$ being arbitrary, we follow the approach of \cite{Breakwelletal1966} and choose $J_{1,1}=J_{1,2}=0$, implying that
\[ K_1(2in\omega_1-G)=\frac{\omega_1K_3}{2}(2n\omega_1-iG), \quad K_2(G-2in\omega_2)=\frac{\omega_2K_4}{2}(2n\omega_2-iG)\]
 i.e.
 \begin{eqnarray}\frac{K_1}{\omega_1(2n\omega_1-iG)}=\frac{K_3}{2(2in\omega_1-G)}=h_1\ \mbox{(say)}\label{eq:k1k3h1}\\
  \frac{K_2}{\omega_2(2n\omega_2-iG)}=\frac{K_4}{2(G-2in\omega_2)}=h_2\ \mbox{(say)}\label{eq:k2k4h2}\end{eqnarray}
 Using equations (~\ref{eq:4ik1}),(~\ref{eq:4ik2}),(~\ref{eq:k1k3h1}) and (~\ref{eq:k2k4h2}) we observe that
 \begin{equation}
 h_j=\frac{1}{2\omega_jM_j\overline{M}_j(M_j^*)^2} \label{eq:hj}
 \end{equation}
 where \begin{eqnarray}M_j&=&(\omega_j^2-2F+n^2)^{1/2},\  M_j^*=(\omega_j^2-2E+n^2)^{1/2},\nonumber\\
\overline{M_j}&=&\sqrt{2} (\omega_j^2-E-2-n^2)^{1/2}, (j=1,2)
\end{eqnarray}
       It is now verified that $H_2$ takes the form:
       \begin{equation} H_2=\frac{1}{2}{(P_1^2-P_2^2+\omega_1^2Q_1^2-\omega_2^2Q_2^2)}\end{equation}
 We observe that
 \begin{equation}
   J=[J_{ij}]_{4\times4}=\left[\begin{array}{c c c c}
  0& 0 &-\frac{M_1}{\omega_1\overline{M_1}}& \frac{i M_2}{\omega_2\overline{M_2}}\\
   -\frac{2n\omega_1}{M_1\overline{M_1}}&  \frac{2n\omega_2}{M_2\overline{M_2}}& \frac{G }{\omega_1M_1\overline{M_1}}& \frac{i G }{\omega_2 M_2\overline{M_2}}\\
  -\frac{\omega_1(m_1^2-2n^2)}{M_1\overline{M_1}}&  \frac{i\omega_2(M_2^2-2n^2)}{M_2\overline{M_2}}& \frac{-nG}{\omega_1 M_1\overline{M_1}}& \frac{n iG }{\omega_2 M_2\overline{M_2}}\\
   \frac{-\omega_1 G }{ M_1\overline{M_1}}& \frac{ i\omega_1 G }{\omega_2 M_2\overline{M_2}}& -\frac{n(\omega_1^2-M_1^2)}{\omega_1 M_1\overline{M_1}}&  -\frac{ni(2\omega_2^2-M_2^2)}{\omega_2 M_2\overline{M_2}} \end{array}\right] \end{equation}
with $l_j^2=4\omega_j^2+9,(j=1,2)$ and $ k_1^2=2\omega_1^2-1,
k_2^2=-2\omega_2^2+1 $.
 Applying a contact transformation from $Q_1,Q_2,P_1,P_2$ to $Q_1',Q_2',P_1',P_2'$ defined by  \cite{Whittaker1965}
\[P_j'=\frac{\partial W}{ \partial Q_j}, \quad Q_j'=\frac{\partial W}{ \partial P_j}, \quad (j=1,2) \]
and
\[W=\sum_{j=1}^2 \left[ Q_j'\sin^{-1}\left( \frac{P_j}{\sqrt{2\omega_jQ_j'}}\right)+\frac{P_j}{2\omega_j}\sqrt{2\omega_j Q_j'-P_j^2}\right]\]
i.e.
\[Q_j=\sqrt{\frac{2Q_j'}{\omega_j}}\cos P_j',  \quad P_j=\sqrt{2\omega_j Q_j'}\sin P_j' \quad (j=1,2)\]
The Hamiltonian $H_2$ is transformed into the form
\[H_2=\omega_1Q_1'-\omega_2Q_2'\]
Denoting the angular variables $P_1'$ and $P_2'$ by $\phi_1$ and $\phi_2$ and the actions $Q_1'$ and $Q_2'$ by $I_1, I_2$ respectively, the second order part of the Hamiltonian transformed into  the normal form
\begin{equation}H_2=\omega_1I_1-\omega_2I_2\label{eq:H2Normal}\end{equation}
The general solution of the corresponding equations of motion are

\begin{equation} I_i={const.}, \quad \phi_i=\pm \omega_i+{const},\  (i=1,2) \label{eq:intmt} \end{equation}

If the oscillations about $L_4$ are exactly linear, the
equation(~\ref{eq:intmt}) represent the integrals of motion and the
corresponding orbits are  given by
\begin{eqnarray}x&=&J_{13}\sqrt{2\omega_1I_1}\cos{\phi_1}+J_{14}\sqrt{2\omega_2I_2}\cos{\phi_2}\label{eq:xb110}\\
y&=&J_{21}\sqrt{\frac{2I_1}{\omega_1}}\sin{\phi_1}+J_{22}\sqrt{\frac{2I_2}{\omega_2}}\sin{\phi_2}\nonumber\\&&+J_{23}\sqrt{2I_1}{\omega_1}\cos{\phi_1}+J_{24}\sqrt{2I_2}{\omega_2}\sin{\phi_2}\label{eq:yb101}
\end{eqnarray}
\section{Conclusion}
We have seen  there are closed curves around the $L_{4(5)}$ so they are stable but the stability range reduced due to radiation effect. The effect of  oblateness and mass of the belt is presented.  We have found the normal form of the second order part of the Hamiltonian. For this we have solved
the aforesaid set of equations. Under the normality conditions, we
have applied a transformation defined by  \cite{Whittaker1965}. We
have also utilized the approach of \cite{Breakwelletal1966} for
reducing the second order part of the Hamiltonian into the normal
form. We  have found that  the second order
part $H_2$ of the Hamiltonian is transformed into the normal form
$H_2=\omega_1I_1-\omega_2I_2.$
\subsection{Acknowledgments}
I would like to acknowledge and extend my heartfelt gratitude to Dr. Uday Dolas for her vital encouragement and support. Especially, I would like to give my special thanks to my wife Mrs. Snehlata Kushwah(Pinky) whose patient love enabled me  to complete this research article.
\section*{References}
\bibliographystyle{jphysicsB} 

\end{document}